\documentclass[11pt]{article}

\usepackage{fullpage,amsmath,amsfonts,amssymb,url}

\newtheorem{lemma}{Lemma}[section]
\newtheorem{theorem}{Theorem}[section]
\newtheorem{proposition}[theorem]{Proposition}

\newcommand{\N}{\mathbb{N}}
\newcommand{\Z}{\mathbb{Z}}

\newcommand{\w}{\mathfrak w}

\newcommand{\legendre}[2]{\genfrac {(}{)}{1pt}{}{#1}{#2}}
\newcommand{\ZNS}[1]{(\Z/{#1}\Z)^{\ast}}

\newcommand{\mat}[4]{\left(\begin{array}{cc} {#1} & {#2} \\ {#3} & {#4}\end{array}\right)}

\begin{document}

\title{Modular equations for some $\eta$-products}
\author{
Fran\c{c}ois Morain \\
INRIA Saclay--\^Ile-de-France \\
\& Laboratoire d'Informatique (CNRS/UMR 7161) \\
\'Ecole polytechnique \\
91128 Palaiseau \\
France \\
morain@lix.polytechnique.fr}
\date{February 8, 2011}

\maketitle

\begin{abstract}
The classical modular equations involve bivariate polynomials that can
be seen to be univariate with coefficients in the modular invariant
$j$. Kiepert found modular equations relating some $\eta$-quotients and
the Weber functions $\gamma_2$ and $\gamma_3$. In the present work, we
extend this idea to double $\eta$-quotients and characterize all the
parameters leading to this kind of equation. We give some properties
of these equations, explain how to compute them and give numerical
examples.
\end{abstract}


\section{Introduction}

Let $\eta$ denote Dedekind's function. When $N>1$ is an integer,
$\eta$-quotients of the form $f = \prod_{d\mid N}\eta(z/d)^{r_d}$ are
functions for $\Gamma^0(N)$ when the integer $r_d$'s satisfy some
properties known as Newman's Lemma \cite{Newman57}. In other words,
there exists a bivariate polynomial $\Phi[f](X, J)$ such that
$\Phi[f](f(z), j(z)) = 0$ for all $z$, where $j$ is the classical
modular invariant. 

In some cases, there exist equations of the form $\Phi[f](X, G_3, G_2)$
where $\Phi[f](f(z), \gamma_3(z), \gamma_2(z)) = 0$ for the Weber function
$\gamma_3$, $\gamma_2$. Kiepert was the first to compute modular
equations of this type for $f = \w_p = \eta(z/p)/\eta(z)$ for $p \leq
29$ (see \cite{Kiepert86}). Weber cites some examples in
\cite[\S 72]{Weber08-III} and Antoniadis \cite{Antoniadis84} extended
this to $p \leq 61$.

In the present work, we study such equations for the
double $\eta$-quotients $\w_{p_1, p_2}^e$, as introduced in
\cite{EnSc04}. We give all parameters $(p_1, p_2, e)$ leading to
equations in $\gamma_2$ and $\gamma_3$.

Section 2 recalls known facts on Weber and $\eta$ functions. Section 3
deals with the case of $\w_p$ where we introduce a faster variant of the
classical algorithm to compute the modular equation via series
expansions. Section 4 proves the necessary results for $\w_{p_1,
p_2}$, gives algorithms to compute the equations in the spirit of
Section 3, and we add numerical examples.

\medskip
\noindent
{\bf Notations:} If $u$ is some function, we will note $\Phi[u](X, J)$ the
corresponding modular equation. If $u= j(nz)$, we will note $\Phi_n$
to simplify.

\section{Preliminaries}

\subsection{Properties of the functions $\gamma_2$ and $\gamma_3$}

We will use the traditional notations
$$T = \mat{1}{1}{0}{1}, \quad S = \mat{0}{-1}{1}{0}$$
and use the notation $f \circ M$ to denote the function $z \mapsto f(M
z)$. The modular invariant is $j(z) = j(q) = 1/q + 744 + \cdots$ with
$q = \exp(2i\pi z)$. The classical Weber functions are
$$\gamma_2(q) = j(q)^{1/3} = q^{-1/3} (1 + 248 q + 4124 q^2 + 34752
q^3 + 213126 q^4 + O(q^5)) \in q^{-1/3} (1 + \Z[[q]]),$$
$$\gamma_3(q) = (j(q)-1728)^{1/2} = q^{-1/2} (1
-492\,q-22590\,{q}^{3}+O \left( {q}^{5} \right)) \in q^{-1/2} (1 +
\Z[[q]]).$$
If $n$ is an integer, we note $\zeta_n = \exp(2 \pi i / n)$.
Remember that $j$ is invariant through $T$ and $S$ and that
\begin{equation}\label{action-gamma2}
\gamma_2 \circ T = \zeta_3^{-1} \gamma_2, \quad \gamma_2 \circ S =
\gamma_2,
\end{equation}
\begin{equation}\label{action-gamma3}
\gamma_3 \circ T = - \gamma_3, \quad \gamma_3 \circ S = -\gamma_3.
\end{equation}
Moreover, we have \cite[\S 55]{Weber08-III}
\begin{theorem}\label{th-fdal}
(a) Any function invariant by $T$ and $S$ is a rational function of
$j$.

(b) Any function $f$ satisfying $f\circ T = -f$ and
$f\circ S = -f$ is equal to $\gamma_3$ times a rational function of $j$.

(c) Any function $f$ satisfying $f\circ T = \zeta_3^{\mp 1} f$ and
$f\circ S = f$ is equal to $\gamma_2^{\pm 1}$ times a rational function of $j$.

(d) Any function $f$ satisfying $f\circ T = -\zeta_3^{\mp 1} f$ and
$f\circ S = -f$ is equal to $\gamma_3\gamma_2^{\pm 1}$ times a rational
function of $j$. (Note that $-\zeta_3^{\mp 1} = \zeta_6^{\pm 1}$.)
\end{theorem}
Let us precise this result in a special case.
\begin{proposition}\label{recog}
Let $\mathcal{T}(q)$ be invariant as in Theorem \ref{th-fdal}. Suppose
that $\mathcal{T}(q) \in q^{-a/b} \Z[[q]]$ for an irreducible fraction
$a/b > 0$ with $b \mid 6$. When $b=1$, $\mathcal{T}(q)$ is a
polynomial in $j(q)$.
When $b=2$, $\mathcal{T}(q) / \gamma_3(q)$ is a
polynomial in $j(q)$. For $b=3$, $\mathcal{T} / \gamma_2^i$ is a
polynomial in $j(q)$ where $i \equiv -a\bmod 3$. For $b=6$,
$\mathcal{T} / (\gamma_2^i \gamma_3)$ is a polynomial in $j(q)$,
where $i \equiv -(a+3)/2 \bmod 3$. In all cases, the polynomial in
$j(q)$ has integer coefficients.
\end{proposition}

\noindent
{\em Proof:} in all cases, the integer $i$ is chosen in such a way
that the resulting series $\mathcal{T}'$ is invariant through $S$ and
$T$, therefore a rational function in $j$. Noting that $\mathcal{T}'$
has integer coefficients, by the Hasse principle, so does the
polynomial. $\Box$

From the algorithmic point of view, we have to recognize a polynomial
with integer coefficients applied to $j(q)$, given the first terms of
the series $\mathcal{T}(q)$. Note that we need the order of this
series to be $>0$. We proceed step by step.

\medskip
\noindent
{\bf function} {\sc RecognizePolyInJ}($\mathcal{T}$)

{\sc Input:} a series $\mathcal{T} = c_{v} q^{v} + \cdots + O(q^1)$ with
integer coefficients, $v \leq 0$ and $c_v \neq 0$.

{\sc Output:} a polynomial $P(X)$ of degree $-v$ such that
$\mathcal{T} = P(j(q))$.

\medskip
1. $\mathcal{R} := \mathcal{T}$; $i:=\mathrm{valuation}(\mathcal{R})$; $P:=0$;

2. while $i \leq 0$ do

\ \ \{ at this point $\mathcal{R} = r_i q^i + \cdots + O(q^1)$ with $r_i
\neq 0$ \}

\ \ \ \ \ \ 2.1 $P:=P + r_i X^{-i}$;

\ \ \ \ \ \ 2.2 $\mathcal{R} := \mathcal{R} - r_i j(q)^{-i}$;

\ \ \ \ \ \ 2.3 $i:=\mathrm{valuation}(\mathcal{R})$;

3. return $P$.

\medskip
Note that we can precompute the powers of $j(q)$ whenever needed, so
that each call to the function requires $O(v^2)$ operations. In
large cases, computations can be done using results computed modulo
small primes and reconstructed via the CRT (as done by Atkin, see
\cite{Morain95a}).

\subsection{Formulas for the $\eta$-function}

The following is taken from \cite{EnSc05} and will be our main tool in
the computations of Section \ref{sct-actions}.
\begin{theorem}
\label{th:transformation}
Let $M = \begin {pmatrix} a & b \\ c & d \end {pmatrix} \in \Gamma$ be
normalised such that $c \geq 0$, and $d >0$ if $c = 0$. Write $c = c_1
2^{\lambda (c)}$ with $c_1$ odd; by convention, $c_1 = \lambda (c) = 1$ if
$c = 0$. Define
\[
\varepsilon (M) = \legendre {a}{c_1} \zeta_{24}^{a b + c (d (1 - a^2) - a)
+ 3 c_1 (a - 1) + \frac {3}{2} \lambda (c) (a^2 - 1)}.
\]
For $K \in \N$ write
\[
u_K a + v_K K c = \delta_K = \gcd (a, Kc) = \gcd (a, K).
\]
Then
\[
\eta \left( \frac {z}{K} \right)
 \circ M = \varepsilon \left(\begin {pmatrix}
\frac {a}{\delta_K} & -v_K \\ \frac {Kc}{\delta_K} & u \end
{pmatrix}\right)
 \sqrt{\delta_K (c z + d)} \, \eta \left( \frac {\delta_K z + (u_K b + v_K K
d)}{\frac {K}{\delta_K}} \right),
\]
where the square root is chosen with positive real part.
\end{theorem}
We can decompose the formula into several parts:
$\varepsilon(M) = \mathrm{Jac}(M)
\zeta_{24}^{\mathcal{E}(M)}$ where we distinguish the Jacobi symbol
part and the exponent of $\zeta_{24}$; then, we have the squareroot
part $\mathcal{Q}(M)$ and the $\eta$-part $\mathcal{N}(M)$. When
dealing with a $\eta$-quotient, the aboves formulas are applied by
multiplicativity on the different pieces $\eta(z/d)$ (see below).

\section{Generalized Weber functions}
\label{sct-weber}

\subsection{Definition and properties}

Let $N>3$ be an odd integer. For all factorizations $N
= a d$, let $e = \gcd(a, d)$ and consider the functions
$$P_{c, d, a} = i^{(a-1)/2} \legendre{c}{e} \sqrt{d}
\frac{\eta((c+dz)/a)}{\eta(z)}$$
for $0 \leq c < a$ with $\gcd(c, e)=1$. These functions were
introduced in \cite[\S 72]{Weber08-III}. It is easy to see that
$P_{0, 1, N} = i^{(N-1)/2} \eta(z/N)/\eta(z) = i^{(N-1)/2} \w_N(z)$
where the function $\w_N$ was studied in \cite{EnMo09}.

Weber proves that in all cases,
$P_{c, d, a}^{24}$ are roots of a modular equation. In some cases, the
results are better, for instance:
\begin{theorem}
If $\gcd(N, 6) = 1, 12 \mid c$, then the
$P_{c, d, a}^2\gamma_2^{N-1}\gamma_3^{(N-1)/2}$ are roots of a modular
equation.
\end{theorem}

\subsection{Computations in the prime order case}
\label{ssct-wp}

For a prime $N = p > 3$, this setting simplifies to
$$x_{0, p, 1} = p \left(\frac{\eta(p z)}{\eta(z)}\right)^2,
\quad x_{12h, 1, p} = (-1)^{(p-1)/2}
\left(\frac{\eta\left(\frac{12 h +z}{p}\right)}{\eta(z)}\right)^2,
 0 \leq h < p.$$

\begin{theorem}
The numbers $x_{c, d, a} \gamma_2(z)^{p-1}
\gamma_3(z)^{(p-1)/2}$ are roots of a modular equation whose
coefficients are rational functions of $j(z)$. In particular, the
constant term is $(-1)^{(p-1)/2} p$.
\end{theorem}

Antoniadis extended the results of Kiepert to
$p \leq 61$ and gave more properties of the polynomials
\cite{Antoniadis84}. He computed
the equation by solving a linear system in the unknown coefficients of
the equation, using the $q$-expansion of $j(q)$ and the fact that
$x_{0, p 1}$ must be a root of the equation.

A standard approach (already known to Enneper \cite[\S 52]{Enneper90}) is to
compute the powers sums of the roots of the equation, recognize them
as polynomials in our variables, and then terminate using the
classical Newton relations.
Inspecting our roots, we see that the $q$-expansion of $x_{0, p, 1}$
has positive order, and all $x_{12h, 1, p}$ have negative order. So
the power sums can be computed using the $x_{12h, 1, p}$ only; we can find
formulas for the $q$-expansion of $\sum_{h=0}^{p-1} x_{12h, 1, p}^k$
if needed.

A better approach is to look at the reciprocal polynomial, whose roots
are the $1/x_{0, p, 1}$ and $1/x_{12h, 1, p}$ and only the first one
contributes to the power sums. Write $(p-1)/12 = e'/\delta$ as an
irreducible fraction with $6\mid \delta$. Noting that
$$p/x_{0, p, 1} = q^{(1-p)/12} (1 + \cdots) = q^{-e'/\delta} (1 + \cdots),$$
we see that all powers are expressible as functions whose expansions
satisfy Proposition \ref{recog}.

\noindent
The algorithm is:

1. compute $S_k = p/x_{0, p, 1}^k$ and recognize it as a polynomial in
the usual variables.

2. use Newton's formulas.

3. Remove the powers of $p$.

Note that the largest power is $(p/x_{0, p, 1})^{p+1} =
q^{-(p^2-1)/12} (1 + \cdots)$ where the exponent is an integer as soon
as $p > 3$. Therefore, we need up to $(p^2-1)/12$ terms in the
$j$-series.


We have
$$S_1 = 11/x_{0, 11, 1} = q^{-5/6} - 2 q^{1/6} - q^{7/6} + 2 q^{13/6}
+ O(q^{19/6}).$$
Dividing by $\gamma_2\gamma_3$, we find
$$1 + 242 q + O(q^2)$$
which must be a polynomial in $j(q)$, hence the constant 1. The other
coefficients are given in Table \ref{wp11}. We have replaced
$\gamma_2$ (resp. $\gamma_3$) by $G_2$ (resp. $G_3$).
\begin{table}[p]
$$\begin{array}{|r|l|}\hline
k & (11/x_{0, 11, 1})^k \\ \hline
2 & q^{-5/3} - 4 q^{-2/3} + 2 q^{1/3} + \cdots = G_2^2 (J - 1244) \\
3 &q^{-5/2} - 6 q^{-3/2} + 9 q^{-1/2} + 10 q^{1/2} + \cdots 
= G_3(J^2 - 1002 J + 59895) \\
4 & q^{-10/3} - 8 q^{-7/3} + 20 q^{-4/3} - 70 q^{2/3} + \cdots
= G_2 (J^3- 2488 J^2 + 1510268 J - 135655520) \\
5 & q^{-25/6} - 10 q^{-19/6} + 35 q^{-13/6} - 30 q^{-7/6} -
105 q^{-1/6} + 238 q^{5/6} + \cdots \\
& = G_3 G_2^2 (J^3 - 2246 J^2 + 1287749 J -145411750) \\
6 & q^{-5} - 12 q^{-4} + 54 q^{-3} - 88 q^{-2} - 99 q^{-1} + 540 - 418
q + \cdots \\
 & = J^5 - 3732 J^4 + 4586706 J^3 - 2059075976 J^2 + 253478654715 J - 
2067305393340 \\
7 & q^{-35/6} - 14 q^{-29/6} + 77 q^{-23/6} - 182 q^{-17/6} +
924 q^{-5/6} - 1547 q^{1/6} + \cdots \\
& = G_3 G_2 (J^5 - 3490 J^4 + 4063139 J^3 -
1796527998 J^2 + 247854700555 J - 4740750382830) \\
8 & q^{-20/3} - 16 q^{-17/3} + 104 q^{-14/3} - 320 q^{-11/3} +
260 q^{-8/3} + 1248 q^{-5/3} - 3712 q^{-2/3} \\
& + 1664 q^{1/3} + \cdots 
 = G_2^2 (J^6 - 4976 J^5 + 9210680 J^4 - 7786404608 J^3 + 2955697453292 J^2 \\
& \hspace*{4.5cm} - 418137392559040 J + 12629117378938720) \\
9 & q^{-15/2} - 18 q^{-13/2} + 135 q^{-11/2} - 510 q^{-9/2} +
765 q^{-7/2} + 1242 q^{-5/2} - 7038 q^{-3/2} + 8280 q^{-1/2} \\
& \quad + 9180 q^{1/2} + \cdots =
G_3 (J^7 - 4734 J^6 + 8386065 J^5 - 6877048710 J^4 +
2611195915626 J^3 \\
& \hspace*{3.5cm}
- 398512009001700 J^2 +16457557949779815 J - 41283301866181650) \\
10 & q^{-25/3} - 20 q^{-22/3} + 170 q^{-19/3} - 760 q^{-16/3} + 
1615 q^{-13/3} + 476 q^{-10/3} - 11210 q^{-7/3} \\
& \quad + 22440 q^{-4/3} +1615 q^{-1/3}- 64600 q^{2/3}+\cdots \\
& = G_2 (J^8 - 6220 J^7 + 15382190 J^6 -
19242776200 J^5 + 12809764457825 J^4 \\
&\quad\quad - 4368737795118764 J^3 + 669619352632925750 J^2 \\
&\quad\quad -33921007872189625000 J + 233702090524237500000) \\
11 & q^{-55/6} - 22 q^{-49/6} + 209 q^{-43/6} - 1078 q^{-37/6} + 
2926 q^{-31/6} - 1672 q^{-25/6} - 15169 q^{-19/6} \\
& + 47234 q^{-13/6} - 31350 q^{-7/6} - 107426 q^{-1/6} + 218680 q^{5/6}+\cdots \\
& G_3 G_2^2 (J^8 - 5978 J^7 + 14256527 J^6 - 17312108670 J^5 +
11327366012605 J^4 \\
& - 3889904574252522 J^3 + 631138185556080950 J^2 \\
& - 38141443583282670180 J + 473098671409604281800) \\
12 & q^{-10} - 24 q^{-9} + 252 q^{-8} - 1472 q^{-7} + 4830 q^{-6} -
6048 q^{-5} - 
16744 q^{-4} + 84480 q^{-3} \\
&  - 113643 q^{-2} - 115920 q^{-1} + 534612 - 370920 q +\cdots \\
& J^{10} - 7464 J^9 + 23101236 J^8 - 38353325536 J^7 +
36913772324730 J^6 - 20784851556729552 J^5 \\
& \quad + 6580486714450069928 J^4 - 1063011399511905159360 J^3
+ 72005127765018136775955 J^2\\
& \quad - 1322204967509387392211000 J + 1424583710586688670191932 \\
\hline
\end{array}$$
\caption{Computations for $p=11$.\label{wp11}}
\end{table}
The corresponding polynomial is (after reductions between variables)
$$\Phi[-\w_{11}^2](F, G_2) = 
F^{12} - G_3 G_2 F^{11} - 242 G_2^2 F^{10} - 19965 G_3 F^9$$
$$\hspace*{5cm} - 585640 G_2 F^8 + 159440490 F^6 - 285311670611.$$
Taking its reciprocal and removing the spurious powers of $p$ yields:
$$F^{12} - 990 F^6 + 440 G_2 F^4 + 165 G_3 F^3 + 22 G_2^2 F^2 + G_3
G_2 F - 11,$$
already computed by Weber.

Note that one drawback of the approach is the large degree and sizes of the
coefficients before reduction via Newton formulas. However, if
computations are performed using CRT primes, this is not a problem,
since we compute the final polynomial modulo the primes.

The smallest cases are
$$\Phi[\w_5^2](X, G_2) = {X}^{6}+10\,{X}^{3} - G_2 X +5,$$
$$\Phi[-\w_7^2](X, G_3) = X^8+14 X^6+63 X^4+70 X^2+G_3 X-7,$$
$$\Phi[\w_{13}^2](X, J) =
{X}^{14}+26\,{X}^{13}+325\,{X}^{12}+2548\,{X}^{11}+13832\,{X}^{10}+
54340\,{X}^{9}+157118\,{X}^{8}+333580\,{X}^{7}$$
$$+509366\,{X}^{6}+534820
\,{X}^{5}+354536\,{X}^{4}+124852\,{X}^{3}+15145\,{X}^{2}+ \left( 746-J
 \right) X+13.$$

\bigskip
\noindent
{\bf Remark.} We concentrated here on the prime index case. The same
work can be done on composite ones. Note also that we could use
resultants for that task, noting the following.
Suppose $p$ is prime and $M$ is an integer prime to $p$; write $N = p
M$. Write
$$\w_{p M}^s(z) = \left(\w_{p}(z) \w_{M}(z/p)\right)^s.$$
On the other hand:
\begin{eqnarray*}
\Phi[\w_{p}^{s_1}](\w_{p}^{s_1}(z), j(z)) &=& 0,\\
\Phi[\w_{M}^{s_2}](\w_{M}^{s_2}(z/p), j(z/p)) &=& 0,\\
\Phi_{p}(j(z), j(z/p)) &=& 0.
\end{eqnarray*}
Writing $Z = \w_{p M}^s(z)$, $X = \w_{p}(z)$, $Y =
\w_{M}(z/p)$, the different quantities are related via the
algebraic equations:
$$Z = X^s Y^s,$$
$$\Phi[\w_{p}^{s_1}](X^{s_1}, J) = 0,$$
$$\Phi[\w_{M}^{s_2}](Y^{s_2}, J') = 0,$$
$$\Phi_{p}(J, J') = 0,$$
and the variables can be eliminated via resultants to get a modular
equation in $Z$ and $J$, that needs to be factored to get the correct
polynomial.

\section{Double $\eta$-quotients}
\label{sct-actions}

\subsection{Definition and statement of the result}

For primes $p_1$ and $p_2$, let
$$\w_{p_1, p_2}^s = \left( \frac {\eta \left( \frac {z}{p_1} \right)
\eta \left( \frac {z}{p_2} \right)}{\eta \left( \frac {z}{p_1 p_2}
\right)
\eta (z)} \right)^s =\left(\frac{\w_{p_1}(z)}{\w_{p_1}(z/p_2)}\right)^s$$
where $s = \frac {24}{\gcd (24, (p_1 - 1)(p_2 - 1))}$ is the smallest
integer such that $s r$ is an integer, where $r = (p_1 - 1)(p_2 -
1)/24$. Note that $s \mid 24$; and $s \mid 6$ when $p_1$ and $p_2$ are odd
primes. It is shown in \cite{EnSc05} that the function $\w_{p_1,
p_2}^s$ is a function on $\Gamma^0(p_1 p_2)$; properties of the
classical modular equation are also given.

We can now state the result that we will prove in this Section.
\begin{theorem}\label{th-fdal2}
Let $p_1, p_2$ be two primes, $N = p_1 p_2$, $s = 24/\gcd (24, (p_1 -
1)(p_2 - 1))$, $e \neq s$ a divisor of $s$ and $\delta = s/e$. If $N\equiv
1\bmod \delta$ and the parameters are chosen in Table
\ref{thp1p2}, then there exists a modular equation
$\Phi[(-1)^{\delta+1} \w_{p_1, p_2}^e]$ whose coefficients are
rational functions in $\gamma_3$, $\gamma_2$.

\begin{table}[hbt]
$$\begin{array}{|c|c||c|c|c|}\hline
p_1 & p_2 & s & e & \delta \\ \hline
2  & 2 & 24 & 8 & 3 \\
2  & 5\bmod 12 & 6 & 2 & 3 \\
2 & 11\bmod 12 & 12 & 4 & 3 \\
\hline
3 & 3 & 6 & 3 & 2 \\
3 & 7\bmod 12 & 2 & 1 & 2 \\
3 & 11\bmod 12 & 6 & 3 & 2 \\ 
\hline
5\bmod 12 & 5\bmod 12 & 3 & 1 & 3 \\ 
5\bmod 12 & 11\bmod 12 & 3 & 1 & 3 \\ 
\hline
7\bmod 12 & 7\bmod 12 & 2 & 1 & 2 \\ 
7\bmod 12 & 11\bmod 12 & 2 & 1 & 2 \\ 
\hline
11\bmod 12 & 11\bmod 12 & 6 & 1 & 6 \\ 
\hline
\end{array}$$
\caption{Values of $p_1$ and $p_2$ leading to a modular equation
$\Phi[(-1)^{\delta+1} \w_{p_1, p_2}^e]$.\label{thp1p2}}
\end{table}
\end{theorem}

The following Lemma is used in the Theorem.
\begin{lemma}\label{lemp}
Let $\delta \in \{2, 3, 6\}$ be as above and suppose $N = p_1
p_2\equiv 1\bmod \delta$. Then $p_i \equiv -1\bmod \delta$.
\end{lemma}

\noindent
{\em Proof:} 
For $\delta = 2$, $N\equiv 1\bmod 2$ gives the
answer. When $3\mid\delta$, we cannot have $p_i=3$ since $N \equiv
1\bmod \delta$. For $\delta$ to be equal to $3$ (resp. $6$), surely we
cannot have $p_i \equiv 1\bmod 3$ (resp. $6$). This leaves $p_i\equiv
-1\bmod 3$ (resp. $6$). $\Box$

The proof of the Theorem will take use several intermediate results
that we will present in as much a compact way as possible.
When $p_1 \neq p_2$, we will make the convention that $p_1$ is
odd (so that we may have $p_2 = 2$). Moreover, we let $u$ and $v$ be
two integers such that $u p_1 + v p_2 = 1$. To simplify the proofs, we
will be mostly looking at properties using $p_2$, this case being
complicated when $p_2 = 2$. Reciprocally, using $p_1$ and $p_2$
supposes that $p_1 \neq p_2$.
The results and proofs are of course symmetrical by
exchanging $p_1$ and $p_2$. In case of equality, we will write $p_1 =
p_2 = p$.

\subsection{The conjugates of $\w_{p_1, p_2}$}

In \cite{EnSc05} are given the conjugates of $\w_{p_1, p_2}^s$ (with
some minor typos). Here, we need precise the expansions of $\w_{p_1,
p_2}$. In view of Theorem \ref{th:transformation}, the value of
$\w_{p_1, p_2} \circ M$ can be composed as
$$\w_{p_1, p_2} \circ M = \mathrm{Jac}(M) \zeta_{24}^{\mathcal{E}(M)}
\mathcal{Q}(M) \mathcal{N}(M)$$
where the first part cumulates Jacobi symbols, the second the
exponents of $\zeta_{24}$, the third one is the product of the
squareroots and the last one the $\eta$ quotient.
To ease notations, we also put $\phi = \zeta_{24}^{24 r} =
\zeta_{24}^{(p_1-1) (p_2-1)}$. We use the notations and philosophy
of computations from \cite{EnSc05}.

\begin{proposition}
Let $p_1$ and $p_2$ be two primes. In all cases, we have the $N+1$
following conjugate functions:
$$\begin{array}{clcc}
M & \multicolumn{1}{c}{\w_{p_1, p_2} \circ M} & \mathrm{ord} & l \\ \hline
T^{\nu} = \mat{1}{\nu}{0}{1} & A_\nu(z) = 
  \frac{
	\eta\left(\frac{z+\nu}{p_1}\right)
	\eta\left(\frac{z+\nu}{p_2}\right)
  }{
	\eta(z+\nu)
	\eta\left(\frac{z+\nu}{p_1 p_2}\right)
  }
& -\frac{r}{p_1 p_2} & \zeta_{N}^{-\nu r} \\
& = \w_{p_1, p_2}(z+\nu), 0 \leq \nu < N & & \\
S = \mat{0}{-1}{1}{0} & B(z) =
  \frac{
        \eta(p_1 z)
        \eta(p_2 z)
  }{
        \eta(z)
        \eta(p_1 p_2 z)
  }
 = \w_{p_1, p_2}(Nz) & - r & 1 \\
\end{array}$$
The remaining $p_1+p_2$ conjugates are:
$$\begin{array}{clcc}
M & \multicolumn{1}{c}{\w_{p_1, p_2} \circ M} & \mathrm{ord} & l \\ \hline
M_{1, \mu} = \mat{\mu p_2}{-1}{1}{0} & C_{1, \nu}(z) = 
\phi^{\theta_1(\nu)} \varepsilon_1
  \frac{
	\eta\left(\frac{z+\nu}{p_1}\right)
	\eta\left(p_2 (z+\nu)\right)
  }{
	\eta(z+\nu)
	\eta\left(\frac{p_2(z+\nu)}{p_1}\right)
  }
, 0 \leq \nu < p_1 & \frac{r}{p_1} &
\phi^{\theta_1(\nu)} \varepsilon_1 \zeta_{p_1}^{\nu r} \\ 
M_{1, 0} = \mat{v p_2}{-u p_1}{1}{1}&&&\\
M_{2, \mu} = \mat{\mu p_1}{-1}{1}{0} & C_{2, \nu}(z) = 
\phi^{\theta_2(\nu)} \varepsilon_2
  \frac{
	\eta\left(p_1(z+\nu)\right)
	\eta\left(\frac{z+\nu}{p_2}\right)
  }{
	\eta(z+\nu)
	\eta\left(\frac{p_1(z+\nu)}{p_2}\right)
  }
, 0 \leq \nu < p_2 & \frac{r}{p_2} &
\phi^{\theta_2(\nu)}\varepsilon_2 \zeta_{p_2}^{\nu r} \\ 
M_{2, 0} = \mat{u p_1}{-v p_2}{1}{1} &&& \\
\end{array}$$
where in the case of $C_{2, \nu}$, we set $\nu \equiv -(\mu p_1)^{-1} \bmod
p_2, v_2 = (1+p_1\mu\nu)/p_2$ for $\mu \neq 0$ (equivalently $\nu \neq
0$; $\mu =0$ corresponds to $\nu = 0$). When $\nu > 0$, we get
$$\theta_2(\nu) =
\left\{
  \begin{array}{ll}
 \mu ((p_2+1) v_2 + 1) + \nu & \text{if } p_2 \neq 2, \\
 (3 p_1 + 2) \frac{\nu+1}{2} & \text{if } p_2 = 2.
  \end{array}
\right.
$$
Moreover
$$\theta_2(0) =
\left\{
  \begin{array}{ll}
u v (p_2+1) + u - 1& \text{if } p_2 \neq 2, \\
\frac{(3 u + 2) (u-1)}{2} & \text{if } p_2 = 2.
  \end{array}
\right.
$$
Also, 
$$\varepsilon_2 = 
\left\{\begin{array}{cl}
\legendre{p_1}{p_2} & \text{if } p_2 \neq 2,\\
1 & \text{if } p_2 = 2.\\
\end{array}\right.$$

When $p_1 = p_2 = p$, we must consider the $p-1$ following conjugate
functions:
$$\begin{array}{clcc}
\text{matrix} & \multicolumn{1}{c}{\w_{p, p}\circ M} & \mathrm{ord} & l \\ \hline
M_\mu = \mat{\mu p}{-1}{1}{0}& C_\nu(z) = 
\sqrt{p} \varepsilon(\nu) \zeta_{24}^{\theta(\nu)}
  \frac{\eta(p z)^2}{\eta(z) \eta\left(z + \frac{\nu}{p}\right)},
  1\leq \nu < p
& \frac{p-1}{12} & \sqrt{p} \varepsilon(\nu) \zeta_{24 p}^{p\theta(\nu)-\nu} \\
\end{array}$$
where $1 = -\mu \nu + v p$, $\varepsilon(\nu) = \legendre{-\nu}{p}$ if $p$ odd
(resp. $1$ when $p = 2$) and
$$\theta(\nu) = 
 \left\{
  \begin{array}{ll}
p\nu\,(1-{\mu}^{2})+ \left( -3\,p+2+v \right) \mu-3+3\,p
& \text{if } p \text{ is odd},\\
0 & \text{if } p=2. \\
  \end{array}
 \right.
$$
\end{proposition}

\medskip
\noindent
{\em Proof:} the cases of the $A_\nu$ matrices and of $B$ are treated
without difficulty, as in \cite{EnSc05}. The value of $\mathcal{Q}(M)$
is $1$, unless we are dealing with the case $p_1 = p_2 = p$. The
computations for the $C$ matrices involve non-zero exponents for
$\zeta_{24}$.

\medskip
\noindent
{\bf Case $p_1 \neq p_2$:}

In the same lines as in \cite{EnSc05}, we first prove the result for
$C_{2, \nu}$ when $\nu > 0$. Iterate over $1 \leq \mu < p_2$ and define
$\nu = - (\mu p_1)^{-1}\bmod p_2 \in \{1, \ldots, p_2-1\}$,
$v_2 = (1 + \mu \nu p_1)/p_2$. Note that $\nu\mapsto \mu$ is an
involution and the corresponding $v_2$'s are equal. Moreover,
iterating over $1 \leq \mu < p_2$ is the same as iterating over $1\leq
\nu < p_2$. We find
$$\mathcal{N}(M_{2, \mu}) = \frac{\eta(p_1 z)\eta((z+\nu)/p_2)}{\eta(z)\eta(p_1
(z+\nu) / p_2)} = \zeta_{24}^{\nu (1 - p_1)}  \frac{\eta(p_1
(z+\nu))\eta((z+\nu)/p_2)}{\eta(z+\nu)\eta(p_1 (z+\nu) / p_2)}.$$

\noindent
(a) Assume first $p_2 \neq 2$. We compute $\mathrm{Jac}(M_{2, \mu}) =
\legendre{p_1}{p_2}$, and the total exponent of $\zeta_{24}$ is
\begin{eqnarray*}
\nu (1 - p_1) + \mathcal{E}(M_{2, \mu}) &=& \left( p_1-1 \right)  \left( p_2\,\nu\,p_1\,\mu^{2}+p_2\,\nu+2\,\mu\,p_2-\mu\,v_2-\nu-\mu \right) \\
&=& (p_1-1) ((p_2-1) (\nu \mu^2 p_1 + \nu + 2\mu) + \mu (1-v_2 +
\mu\nu p_1))\\
&=& (p_1-1)(p_2-1) (\nu \mu^2 p_1 + \nu + 2\mu + \mu v_2) \\
&=& (p_1-1)(p_2-1) (\mu ((p_2+1) v_2 + 1) + \nu)
\end{eqnarray*}
where we have used $p_2 v_2 = 1 + \mu\nu p_1$ twice.

When $p_2 = 2$, we find $\mathrm{Jac}(M_{2, \mu}) = 1$ and the total exponent
of $\zeta_{24}$ is
$$\nu (1 - p_1) + \mathcal{E}(M_{2, \mu}) = \left(p_1-1 \right)
(3 p_1 \mu^2 (\nu+1) + \mu (3\mu -1) + \nu)/2.$$
Since $\nu$ is odd, $\mu = 1$ and the exponent reduces to
$$(p_1-1) (3 p_1 + 2) \frac{\nu+1}{2}.$$

\medskip
\noindent
(b) For $C_{2, 0}$,
$$\mathcal{N}(M_{2, 0}) = \frac{\eta(p_1(z+1))
\eta(z/p_2)}{\eta(z+1)\eta(p_1 z / p_2)} = \zeta_{24}^{p_1-1}
\frac{\eta(z/p_2) \eta(p_1 z)}{\eta(z) \eta(p_1 z / p_2)},$$

Assume first $p_2 \neq 2$. Then
$\mathrm{Jac}(M_{2, 0}) = \legendre{p_1}{p_2}$ and
the exponent of $\zeta_{24}$ is
\begin{eqnarray*}
p_1-1 + \mathcal{E}(M_{2, 0}) &=& -(p_1-1) ((p_2-1) (p_1 u^2 - 2 u + 1) + u
(p_1 u + v - 1)) \\
&=& - 24 r (p_1 u^2 - 2 u + 1 - u v) \\
&=& 24 r (u v (p_2+1) + u - 1).
\end{eqnarray*}

When $p_2 = 2$, we find $\mathrm{Jac}(M_{2, 0}) = 1$ and the total exponent
of $\zeta_{24}$ is
$$(p_1-1) \frac{(3 u + 2) (u-1)}{2}.$$

\medskip
\noindent
{\bf Case $p_1 = p_2 = p$:}

In all cases:
$$\mathcal{N}(M_\mu) = \sqrt{p} \frac{\eta(p z)^2}{\eta(z) \eta(z + \nu/p)}$$
where $1 = -\mu \nu + v p$.

When $p \neq 2$, we find $\mathrm{Jac}(M_\mu) = \legendre{\mu}{p}$ and
the exponent given by $\theta(\nu)$. When $p = 2$, $\mathrm{Jac}(M_\mu) =
1$ and the exponent given by $v-1 = 0$. $\Box$

\subsection{Action of $T$ and $S$}\label{ssct:g23action}

This section is devoted to the proofs of the actions of $T$ and $S$ on
our basic functions as stated in the following two propositions.

\begin{proposition}\label{T-action}
(i) $B \circ T = \phi^{-1} B$.

(ii) For $0 \leq \nu < N-1$, we have $A_\nu \circ T = A_{\nu+1}$;
$A_{N-1} \circ T = \phi^{-1} A_0$.

(iii) For $0 \leq \nu < p_2-1$, $C_{2, \nu} \circ T =
\phi^{\theta_2(\nu)-\theta_2(\nu+1)} C_{2, \nu+1}$; $C_{2, p_2-1}
\circ T = \phi^{\theta_2(p_2-1)-\theta_2(0) + 1} C_{2, 0}.$

(iv) For $1 \leq \nu < p$, $C_{\nu} \circ T = \zeta_{24}^{2p-2} C_{\nu}$.
\end{proposition}

\medskip
\noindent
{\em Proof of Proposition \ref{T-action}:}

\noindent
\iftrue 
(i), (ii) and (iv) are direct applications of Theorem
\ref{th:transformation}.
\else
(i) we find
$$B \circ T = \zeta_{24}^{p_1 + p_2 - 1 - p_1 p_2} B =
\phi^{-1} B.$$

\medskip
\noindent
(ii) For $0 \leq \nu < N-1$, we have $A_\nu \circ T = A_{\nu+1}$. For $\nu
= N-1$, we find
$$A_{N-1} \circ T =   \frac{
	\eta\left(\frac{z+N}{p_1}\right)
	\eta\left(\frac{z+N}{p_2}\right)
  }{
	\eta(z+N)
	\eta\left(\frac{z+N}{p_1 p_2}\right)
  }
=    \frac{
	\eta\left(\frac{z}{p_1} + p_2\right)
	\eta\left(\frac{z}{p_2} + p_1\right)
  }{
	\eta(z+N)
	\eta\left(\frac{z}{p_1 p_2} + 1\right)
  }
= \zeta_{24}^{p_1+p_2-N-1} A_0 = \phi^{-1} A_0.$$
\fi

\medskip
\noindent
(iii) For $0 \leq \nu < p_2-1$, one has 
$C_{2, \nu} \circ T = \phi^{\theta_2(\nu)-\theta_2(\nu+1)} C_{2, \nu+1}$.
For $\nu = p_2-1$:
$$C_{2, p_2-1} \circ T = \phi^{\theta_2(p_2-1)} \varepsilon_2
  \frac{
	\eta\left(\frac{z+p_2}{p_2}\right)
	\eta\left(p_1 (z+p_2)\right)
  }{
	\eta(z+p_2)
	\eta\left(\frac{p_1(z+p_2)}{p_2}\right)
  }
= \phi^{\theta_2(p_2-1)} \varepsilon_2
  \frac{
	\eta\left(\frac{z}{p_2} + 1\right)
	\eta\left(p_1 z + N\right)
  }{
	\eta(z+p_2)
	\eta\left(\frac{p_1 z}{p_2} + p_1\right)
  }$$
$$= \phi^{\theta_2(p_2-1)} \varepsilon_2 \zeta_{24}^{1+N-p_1-p_2}
  \frac{
	\eta\left(\frac{z}{p_2}\right)
	\eta\left(p_1 z\right)
  }{
	\eta(z)
	\eta\left(\frac{p_1 z}{p_2}\right)
  }
 = \phi^{1 + \theta_2(p_2-1)-\theta_2(0)} C_{2, 0}. \Box$$

\begin{proposition}\label{S-action}
For all primes $p_1$ and $p_2$, one has

(i) $(A_0, B) \circ S = (B, A_0)$.

(ii) When $0 < \nu < p_1 p_2$ and $\gcd(\nu, p_1 p_2) = 1$:
$$A_\nu \circ S = \phi^{\theta_3(\nu)} A_\omega$$
where $1 = -\omega \nu + v_{12} (p_1 p_2)$ and
$$\theta_3(\nu) = \left\{
 \begin{array}{ll} 
-\omega\nu^2-2\nu+\omega+3 + \nu v_{12} & \text{if } p_2 \neq 2,\\
\omega + \nu (v_{12} (1-2p_1) + 2) & \text{if } p_2 = 2. \\
 \end{array}\right.
$$

\noindent
Suppose from now on that $p_1 \neq p_2$. The following hold:

(iii) When $0 < \nu = p_1 \rho < p_1 p_2$:
$$A_\nu \circ S = \phi^{\theta_4(\rho)} C_{2, \varpi}$$
where $1 = -\varpi \nu + w p_2$, and
$$\theta_4(\rho) = \left\{\begin{array}{ll}
- \theta_2(\varpi) + \rho (w (p_2 + 1) + 1) + \varpi & \text{if } p_2
\neq 2, \\
- \theta_2(\varpi) + 3 \frac{p_1+1}{2} \rho^2 + \rho (3 w - 2) + \varpi & \text{if } p_2 = 2.
\end{array}\right.
$$

(iv) When $0 < \nu < p_2$, use $\mu \equiv -1/(\nu p_1) \bmod p_2$ and
$$C_{2, \nu}\circ S = A_{\mu p_1}.$$

(v) We have
$$C_{2, 0} \circ S = C_{1, 0} \times
\left\{\begin{array}{ll}
\phi^{-2 \theta_1(0)} & \text{if } p_2 \neq 2, \\
\phi^{-\theta_1(0) + \left(u^2 p_1 \frac{p_1+1}{2}
+ \frac{1-u}{2}\right)} & \text{if } p_2 = 2. \\
\end{array}\right.
$$
\end{proposition}

\medskip
\noindent
{\em Proof:}

\noindent
(i) We first get:
$$w_{p_1, p_2} \circ \left(T^{\nu} \circ S\right) = w_{p_1, p_2} \circ
\mat{\nu}{-1}{1}{0}$$
and the case $\nu = 0$ yields immediately $A_0 \circ S = B$. On the
other hand, we also have the reassuring result that
$$B \circ S = \w_{p_1, p_2}(-N/z) = \w_{p_1, p_2}(z) = A_0(z).$$

\medskip
\noindent
(ii) When $\gcd(\nu, p_1 p_2) = 1$, we write $1 = -\omega \nu + v_{12} (p_1
p_2)$, and find
$$\mathcal{N}(M_{2, \mu} \circ S) = \frac{\eta((z+\omega)/p_1)
\eta((z+\omega)/p_2)}{\eta(z) \eta((z+\omega)/p_1/p_2)} =
\zeta_{24}^{\omega} A_\omega.$$
When $p_2\neq 2$, $\mathrm{Jac}(M_{2, \mu} \circ S) = 1$ and the total exponent
of $\zeta_{24}$ is
\begin{eqnarray*}
\omega+\mathcal{E}(M_{2, \mu} \circ S) & = & -\omega\nu^2 (24 r - 1)
+\nu(-48 r + 1 + v_{12} (1-p_1-p_2))
+ 24 r (\omega + 3) \\
&=& 24 r (-\omega\nu^2-2\nu+\omega+3) + \nu (\omega\nu+1+v_{12}
(1-p_1-p_2))\\
&=& 24 r (-\omega\nu^2-2\nu+\omega+3 + \nu v_{12})\\
&=& 24 r \theta_3(\nu).
\end{eqnarray*}
When $p_2 = 2$, we also have $\mathrm{Jac}(M_{2, \mu} \circ S) = 1$ and the
exponent of $\zeta_{24}$ becomes
$$(p_1-1) (\omega + \nu (v_{12} (1-2p_1) + 2)).$$
The same type of computations show that the results also holds for
$p_1 = p_2 = 2$.

\medskip
\noindent
(iii) Suppose now that $\nu = \rho p_1$, $1\leq \rho < p_2$. We write
$1 = -\varpi \nu + w p_2$. In all cases
\begin{eqnarray*}
\mathcal{N}(M_{2, \mu} \circ S) &=& \frac{\eta(p_1
z)\eta((z+\varpi)/p_2)}{\eta(z)\eta((p_1 (z+ \varpi))/p_2)} \\
&=& \zeta_{24}^{-p_1\varpi+\varpi}
\frac{\eta(p_1
(z+\varpi))\eta((z+\varpi)/p_2)}{\eta(z+\varpi)\eta((p_1 (z+
\varpi))/p_2)} \\
&=& \zeta_{24}^{-p_1\varpi+\varpi}\left( \phi^{-\theta_2(\varpi)}
\varepsilon_2 C_{2, \varpi}\right).
\end{eqnarray*}
Assume $p_2 \neq 2$. We get $\mathrm{Jac}(M_{2, \mu} \circ S) =
\legendre{p_1}{p_2}$. The partial exponent is given by
\begin{eqnarray*}
-\varpi (p_1-1)+\mathcal{E}(M_{2, \mu} \circ S)
 &=&(p_1-1) (p_2 p_1 \varpi \rho^2 + (-w + 2 p_2-1)\rho + (p_2-1) \varpi)\\
&=& (p_1-1) (p_2 \rho (p_1 \rho \varpi) + (-w + 2 p_2-1)\rho + (p_2-1) \varpi)\\
&=& (p_1-1) (p_2 \rho (w p_2 - 1) + (-w + 2 p_2-1)\rho + (p_2-1) \varpi)\\
&=& (p_1-1) (p_2-1) (\rho (w (p_2 + 1) + 1) + \varpi )
\end{eqnarray*}
yielding the final result.

When $p_2 = 2$, we find $\mathrm{Jac}(M_{2, \mu} \circ S) = 1$ and
$$-\varpi (p_1-1)+\mathcal{E}(M_{2, \mu} \circ S)
 =(p_1-1) \left(\rho (3 w - 2) + 3 \frac{p_1+1}{2} \rho^2 + \varpi\right).$$

\medskip
\noindent
(iv) For $1\leq \nu < p_2$, we compute $\mu \equiv -1/(\nu p_1) \bmod
p_2$ and
$$C_{2, \nu} \circ S = \w_{p_1, p_2} \circ \mat{\mu p_1}{-1}{1}{0} \circ S =
\frac{\eta((z+p_1 \mu)/p_1)\eta((z+p_1 \mu)/p_2)}{\eta(z+p_1
\mu)\eta((z+p_1 \mu)/p_1/p_2)} = A_{p_1 \mu}.$$

\medskip
\noindent
(v) In all cases, we compute
$$\mathrm{Jac}(M \circ S) \mathcal{N}(M \circ S) = 
\legendre{p_2}{p_1} \frac{\eta(z/p_1) \eta(p_2 z-p_2)}{\eta(z-1)
\eta(p_2 z/p_1)} = \legendre{p_2}{p_1} \zeta_{24}^{1-p_2}
\frac{\eta(z/p_1) \eta(p_2 z)}{\eta(z) \eta(p_2 z/p_1)}.$$

When $p_2 \neq 2$, this yields
$$\mathrm{Jac}(M \circ S) \mathcal{N}(M \circ S)= \zeta_{24}^{1-p_2}
\phi^{-\theta_1(0)} C_{1, 0}.$$
The exponent of $\zeta_{24}$ is
\begin{eqnarray*}
1-p_2+ \mathcal{E}(M\circ S)
&=& (p_2-1) (p_1 p_2 v^2 + (u+1-2 p_1) v + p_1-1) \\
&=& (p_1-1) (p_2-1) (-u v (p_1+1) - v + 1) \\
&=& 24 r (-u v (p_1+1) - v + 1) \\
&=& -24 r \theta_1(0)
\end{eqnarray*}
so that $C_{2, 0} \circ S = \phi^{-2 \theta_1(0)} C_{1, 0}$.

When $p_2 = 2$, the exponent of $\zeta_{24}$ becomes
$$-1 + \mathcal{E}(M\circ S) = (p_1-1) \left(u^2 p_1 \frac{p_1+1}{2} +
\frac{1-u}{2}\right),$$
so that the final answer is
$$\phi^{-\theta_1(0) + \left(u^2 p_1 \frac{p_1+1}{2}
+ \frac{1-u}{2}\right)} C_{1, 0}. \quad \Box$$

\medskip
\noindent
\begin{proposition}\label{S-actionp}
We suppose that $p_1 = p_2 = p$. Then

(i) When $\nu = \rho p$, $1 \leq \rho < p$, set $1 = -\varpi \rho + w
p$. Then $A_\nu\circ S = C_\varpi$.

(ii) For all $p$, and all $1 \leq \nu < p$, one has $C_{\nu} \circ S =
A_{\mu p}$ where $\mu \equiv -1/\nu \bmod p$.
\end{proposition}

\medskip
\noindent
{\em Proof:}

\medskip
\noindent
(i) When $p \neq 2$:
\begin{eqnarray*}
A_{\nu} \circ S &=& \sqrt{p} \legendre{\rho}{p} 
\frac{\eta(p z)^2}{\eta(z) \eta(z+\varpi/p)}
\zeta_{24}^{-{\rho}^{2}p\varpi+ \left( -3\,p+2+w \right) \rho+p\varpi-3+3\,p} \\
&=&\sqrt{p} \legendre{-\varpi}{p}
\left(1/\sqrt{p} \varepsilon(\varpi) \zeta_{24}^{-\theta(\varpi)}
C_\varpi\right)
\zeta_{24}^{-{\rho}^{2}p\varpi+ \left( -3\,p+2+w \right) \rho+p\varpi-3+3\,p}\\
&=&\zeta_{24}^{-{\rho}^{2}p\varpi+ \left( -3\,p+2+w \right) \rho+p\varpi-3+3\,p-\theta(\varpi)} C_\varpi \\
&=& C_{\varpi}
\end{eqnarray*}
using
$\theta(\varpi) = p\varpi\,(1-{\rho}^{2})+ \left( -3\,p+2+w \right) \rho-3+3\,p$.

When $p = 2$, $\rho = 1$ implying $\varpi = w = 1$ and
$$A_2 \circ S = \sqrt{2} \zeta_{24}^{w-1} \frac{\eta(2
z)^2}{\eta(z) \eta(z+\varpi/2)} = \zeta_{24}^{-\theta(2)} C_1 = C_1.$$

\medskip
\noindent
(ii) In all cases, we get
$$C_{\nu} \circ S = \frac{\eta((z + p \mu)/p)^2}{\eta(z+p \mu)\eta((z
+ p \mu)/p^2)} = A_{\mu p}. \quad \Box$$

\subsection{Finding invariant functions}

The idea is simple. Using the explicit actions given above, we need to
find suitable modifications of the functions $B$, $A_{\nu}$, $C_{1,
\nu}$, $C_{2, \nu}$ such that the action of $T$ and $S$ on any power sum
coincides with the action on $\gamma_2$, $\gamma_3$ or the product
$\gamma_2\gamma_3$, as given in Section \ref{ssct:g23action}.

Note that $B^e \circ T = \zeta_{24}^{-24 r e} B^e$.
Write $r e = t/\delta$ and remark that this
fraction is irreducible ($s$ being prime to $t$ implies $\delta$ is).
This leads to set $\chi = \phi^{-e} = \zeta_{24}^{-24 r e} =
\zeta_{\delta}^{-t}$, a primitive $\delta$-th root of unity.

The aim of this Section is to prove the following Theorem from which
Theorem \ref{th-fdal2} will follow.
\begin{theorem}\label{th:Sk}
Assume we are in the conditions of Theorem \ref{th-fdal2}.
Define the functions
$$A_\nu' = \chi^{\alpha_0-\nu} A_\nu^e,
\quad B' = B^e;
\quad C_{1, \nu}' = \chi^{\theta_1(\nu)-\nu} C_{1, \nu}^e,
\quad C_{2, \nu}' = \chi^{\theta_2(\nu)-\nu} C_{2, \nu}^e;
\quad C_{\nu}' = \chi^{\mu} C_{\nu}^e,$$
where $\mu \equiv -1/\nu \bmod p$ and
$$\alpha_0 = \left\{\begin{array}{cl}
1 & \text{if } \delta = 2,\\
0 & \text{if } \delta = 3,\\
3 & \text{if } \delta = 6,\\
\end{array}\right.$$
making $\chi^{\alpha_0} = (-1)^{\delta+1} = \chi^{-\alpha_0} = \chi^{-3}$.

Then, for all integers $k$, the quantity
$$\mathcal{S}_k = {B'}^k + \sum_{\nu=0}^{N-1} {A_\nu'}^k + \sum_{\nu=0}^{p_1-1}
{C_{1, \nu}'}^k + \sum_{\nu=0}^{p_2-1} {C_{2, \nu}'}^k = {B'}^k +
\mathcal{S}_{A, k} + \mathcal{S}_{C_1, k} + \mathcal{S}_{C_2, k}$$
satisfies
$\mathcal{S}_k \circ (T, S) = (\chi^k, \chi^{\alpha_0 k})
\mathcal{S}_k$.
\end{theorem}

With these notations, we have
\begin{proposition}\label{prop-fT}
The following hold:

(a) $B' \circ T = \chi B'$.

(b) $\{A_\nu'\}_{\nu} \circ T = \{\chi A_\nu'\}_\nu$.

(c) $\{C_{i, \nu}'\}_\nu \circ T = \{\chi C_{i, \nu}'\}_\nu$.

(d) For all $\nu$, $C_\nu' \circ T = \chi C_\nu'$.
\end{proposition}

\medskip
\noindent
{\em Proof:} (a) and (c) follow easily from Proposition \ref{T-action}.

\noindent
(b) We first obtain $A_{N-1}^e\circ T = \chi A_0^e$. Let us explain
how the choice $A_\nu'$ comes from.
For some function $\alpha$ to be precised later, let us put $A_\nu' =
\chi^{\alpha(\nu)} A_\nu^e$, for which
$$A_\nu' \circ T = \chi^{\alpha(\nu)} A_{\nu+1}^e =
\chi^{\alpha(\nu)-\alpha(\nu+1)} A_{\nu+1}',$$
$$A_{N-1}' \circ T = \chi^{\alpha(N-1)} \chi A_0^e =
\chi^{\alpha(N-1)+1-\alpha(0)} A_0'.$$
We must find $\alpha$ s.t.
$$\alpha(\nu) - \alpha(\nu+1) \equiv 1 \bmod \delta, 0 \leq \nu
< N-1,$$
and
$$\alpha(N-1)-\alpha(0)+1 \equiv 1 \bmod \delta.$$
The first set of equations gives us $\alpha(\nu) \equiv \alpha(0) - \nu \bmod
\delta$ and the second $\alpha(0) - (N-1) \equiv \alpha(0) \bmod
\delta$, which is possible only when $N \equiv 1\bmod \delta$. Setting
$\alpha_0 = \alpha(0)$ yields the result. 

\noindent
(d) Proposition \ref{prop-fT} gives us
$C_\nu' \circ T = \zeta_{24}^{2 e(p-1)} C_\nu'$. A glance at Table
\ref{thp1p2} shows that $p^2-1 \equiv 0 \bmod (24/e)$, which implies
$2 e (p-1) \equiv -(p-1)^2 e \bmod 24$ and therefore
$\zeta_{24}^{2 e(p-1)} = \chi$. $\Box$

The actual value of $\alpha$ is in fact dictated by the other
invariance properties that follow.

\medskip
\noindent
{\bf Remark.} This proposition shows at the same time that we cannot
expect some nice $T$-action when $N \not\equiv 1\bmod \delta$.

\medskip
Let us turn our attention to the $S$-action on our candidate
functions, using the notations of Proposition \ref{S-action}.

\begin{proposition}\label{prop-fS}
(i) $(B', A_0') \circ S = \chi^{\alpha_0} (A_0', B')$.

(ii) When $\gcd(\nu, p_1 p_2) = 1$,
$A_\nu' \circ S = \chi^{\alpha_0} A_\omega'$.

(iii) For $\nu = p_1 \rho$,
$A_\nu' \circ S = \chi^{\alpha_0} C_{2, \varpi}'$.

(iv) For $1 \leq \nu < p_2$, $\mu \equiv -1 / (\nu p_1) \bmod p_2$ and
$C_{2, \nu}' \circ S = \chi^{\alpha_0} A_{\mu p_1}'$.

(v) $C_{2, 0}' \circ S = \chi^{\alpha_0} C_{1, 0}'$.

(vi) For $\nu = \rho p$, $1 \leq \rho < p$, set $1 = -\varpi \rho
+ w p$. For all $p$,
$A_\nu' \circ S = \chi^{\alpha_0} C_\varpi'$.

(vii) For $1\leq \nu < p$, setting $\mu \equiv -1/\nu \bmod p$, we have
$C_\nu' \circ S = \chi^{\alpha_0} A_{\mu p}'$.
\end{proposition}

\noindent
{\em Proof:} 

\noindent
(i) We have $B' \circ S = \chi^{-\alpha_0} A_0'$; $A_0' \circ S =
\chi^{\alpha_0} B'$ and the result follows from $\chi^{-\alpha_0} =
\chi{\alpha_0}$.

\medskip
\noindent
(ii) Proposition \ref{S-action} can be rewritten
$$A_\nu' \circ S = \chi^{\omega - \nu - \theta_3(\nu)} A_\omega'$$
and we simplify the exponent using
$1 = -\omega\nu + v_{12} \bmod \delta$, which leads to:
$$
A_\nu' \circ S = A_\omega' \left\{\begin{array}{ll}
\chi^{-3} & \text{if } p_2\neq 2,\\
\chi^{\nu (-3 + v_{12} (2 p_1 - 1))} & \text{if } p_2 = 2.\\
\end{array}\right.
$$
When $p_2\neq 2$, we use $\chi^{-3} = \chi^{\alpha_0}$. The $p_2 = 2$
can occur only for $\delta = 3$, in which case $p_1 \equiv -1\bmod 3$
and the exponent is $\chi$ is $0$.

\medskip
\noindent
(iii) For $\nu = \rho p_1$: we use $1 = -\varpi \nu + w p_2$ to get
$$A_\nu^e \circ S = \chi^{-\theta_4(\rho)} C_{2, \varpi}^e$$
or
$$\chi^{-\alpha_0 + \nu} A_\nu' \circ S =
\chi^{-\theta_4(\rho)} \chi^{\varpi-\theta_2(\varpi)} C_{2, \varpi}'$$
and we need simplify:
$$-\theta_4(\rho)-\theta_2(\varpi)+\alpha_0+\varpi-\nu.$$
Using the definition of $\theta_4$, we get
$$A_\nu' \circ S = C_{2, \varpi}'
\left\{\begin{array}{ll}
 \chi^{\alpha_0 - \rho ((p_2+1) w + p_1 + 1)} & \text{if } p_2 \neq 2,\\
 \chi^{\alpha_0 - \rho \left((3 \frac{p_1+1}{2}) \rho +p_1 + 3 w
 -2\right)} & \text{if } p_2 = 2.\\
\end{array}\right.
$$
and we conclude using $p_i \equiv -1\bmod \delta$.

\medskip
\noindent
(iv) For $1\leq \nu < p_2$, we compute $\mu \equiv -1/(\nu p_1) \bmod
p_2$ and
$$C_{2, \nu}' \circ S = \chi^{\theta_2(\nu)-\nu} A_{\mu p_1}^e
= \chi^{\theta_2(\nu)-\nu-\alpha_0 + \mu p_1} A_{\mu p_1}'.$$
Simplifying the exponent gives 
$$C_{2, \nu}' \circ S = A_{\mu p_1}'
\left\{\begin{array}{ll}
\chi^{-\alpha_0 + \mu ((p_2+1) v_2 + p_1+1)} & \text{if } p_2 \neq 2,\\
\chi^{-\alpha_0 + 4 p_1 + 1} & \text{if } p_2 = 2.\\
\end{array}\right.
$$
where for $p_2 = 2$, we used $\nu = \mu = 1$. We conclude as in (iii).

\medskip
\noindent
(v) When $p_2 \neq 2$, we start from
$$C_{2, 0}^e \circ S = \chi^{2 \theta_1(0)} C_{1, 0}^e$$
from which
$$\chi^{-\theta_2(0)} C_{2, 0}' \circ S = \chi^{\theta_1(0)} C_{1, 0}'$$
or
$$C_{2, 0}' \circ S = \chi^{\theta_1(0) + \theta_2(0)} C_{1, 0}'$$
and the exponent is 
$$u v (p_1 + p_2 + 2) + u + v - 2.$$
This quantity  is $\equiv u+v-2\bmod\delta$
since $p_i\equiv -1\bmod\delta$. Moreover $1 \equiv p_1
(u+v)\bmod\delta$ and finally the exponent is $-3\bmod \delta$.

When $p_2 = 2$, we have
$$\chi^{-\theta_2(0)} C_{2, 0}' \circ S = 
\chi^{-\left(u^2 p_1 \frac{p_1+1}{2} + \frac{1-u}{2}\right)} C_{1, 0}'$$
or
$$C_{2, 0}' \circ S = \chi^{-\frac{(p_1^2+p_1-3) u^2 - 3}{2}} C_{1, 0}',$$
and this is $\chi^{0}$ since this can only happen when $\delta = 3$.

\medskip
\noindent
(vi) Since $1 = -\varpi \rho + w p$, we can write
$$A_\nu' \circ S
= \chi^{\alpha_0 - \nu - \kappa(\varpi)} C_\varpi'.$$
and the result comes from the definition of $\kappa$.

\medskip
\noindent
(vii) $C_\nu^e \circ S = A_{\mu p}^e$ or
$$C_\nu' \circ S = \chi^{\kappa(\nu)} \chi^{-\alpha_0 + \mu p}
A_{\mu p}',$$
and we conclude as in (vi). $\Box$

\subsection{Properties of the modular equation}

From the preceding sections, we see that
$$\Phi(F) = (F-B') \prod_{\nu=0}^{N-1} (F-A_\nu')
\prod_{\nu=0}^{p_1-1} (F-C_{1, \nu}') \prod_{\nu=0}^{p_2-1} (F-C_{2,
\nu}')$$
is a modular equation whose coefficients can be expressed in terms of
$j$, $\gamma_2$ or $\gamma_3$ depending on the value of
$\delta$. Before doing this, we may express these coefficients as
Puiseux series.

\begin{proposition}
With the usual notations:

(a) the coefficient of smallest order of $\Phi$ is $q^{-2 r e}$;

(b) the trace has order $r e$;

(c) when $p_1 \neq p_2$, $\Phi(0) = 1$;

(d) when $p_1 = p_2 = p$, $\Phi(0) = (-p)^{e (p-1)/2}$ when $p$ is
odd and $2^4$ when $p = 2$.
\end{proposition}

\noindent
{\em Proof:} 
(a) the coefficient of smallest order comes from the coefficient of
$F^{\psi(N)-N-1}$ which has the order of $B'
\prod_{\nu=0}^{N-1} A_\nu'$ or
$$-r e + \sum_{\nu=0}^{N-1} \frac{-r e}{N} = -2 r e.$$
When $p_1\neq p_2$, $\psi(N)-N-1 = p_1 + p_2$; when $p_1 = p_2 = p$,
this is $p-1$. Note that all other terms have orders stricly less
than this bound. 

As an example, when $s = e$, the degree of the equation in $J$ is $2 r
s$ and the corresponding term is $J^{2 r s} F^{\psi(N)-N-1}$.

Moreover
$$B' \prod_{\nu=0}^{N-1} A_\nu' = \left(\prod_{\nu = 0}^{N-1}
\chi^{\alpha_0 - \nu} \right) \zeta_N^{-r e N (N-1)/2} q^{-2 r e} (1 +
\cdots)$$
$$= \chi^{N \alpha_0 - N (N-1)/2} \zeta_{24N}^{-24 r e N (N-1)/2}
q^{-2 r e} (1 + \cdots)$$
$$= \chi^{\alpha_0 - N (N-1)/2} \zeta_{24N}^{-24 r e N (N-1)/2} q^{-2
r e} (1 + \cdots),$$
using $N\equiv 1\bmod \delta$. When $N$ is odd, this reduces to
$$\chi^{\alpha_0 - N (N-1)/2 + (N-1)/2} q^{-2 r e} (1 + \cdots)
= \chi^{\alpha_0 - (N-1)^2/2} q^{-2 r e} (1 + \cdots).$$

\medskip
\noindent
(b) The dominant term in the sum of the conjugates is that of $B'$,
namely $q^{-r e}$.

\medskip
\noindent
(c) For $p_1 \neq p_2$:
$$\prod_{\nu=0}^{p_2-1} C_{2, \nu}' = \prod_{\nu=0}^{p_2-1}
\chi^{\theta_2(\nu)-\nu} \chi^{-\theta_2(\nu)}
\varepsilon_2^e \zeta_{p_2}^{\nu r e} q^{e r / p_2} (1 + \cdots)
= \chi^{-p_2 (p_2-1)/2} \varepsilon_2^{p_2 e} \zeta_{p_2}^{r e p_2
(p_2-1)/2} q^{e r} (1 + \cdots).$$
Multiplying all together, we find the norm to be of valuation $0$,
hence a constant
$$\vartheta = \chi^{\alpha_0 - N (N-1)/2} \zeta_{24N}^{-24 r e N (N-1)/2}
\chi^{-p_1 (p_1-1)/2} \varepsilon_1^{p_1 e} \zeta_{p_1}^{r e p_1 (p_1-1)/2}
\chi^{-p_2 (p_2-1)/2} \varepsilon_2^{p_2 e} \zeta_{p_2}^{r e p_2
(p_2-1)/2}$$
$$=(\varepsilon_1^{p_1}\varepsilon_2^{p_2})^e 
\chi^{\alpha_0 - N (N-1)/2-p_1 (p_1-1)/2-p_2 (p_2-1)/2}
\zeta_{24N}^{-24 r e N ( (N-1)/2 - (p_1+p_2-2)/2)}.$$

When $p_2 = 2$ (with $p_1$ odd), we have $\delta = 3$ always,
meaning $\alpha_0 = 0$ and $N\equiv 1\bmod 3$. Therefore, noting that
$e$ is always even:
$$\vartheta = \legendre{2}{p_1}^{p_1 e} \chi^1
\zeta_{24 N}^{-24 r e N (p_1-1)/2} = \chi^{(p_1+1)/2} = 1$$
since $p_1 \equiv -1\bmod 3$.

When $p_2 \neq 2$, both $p_i$ being odd, we may use the quadratic
reciprocity law to find
$$\vartheta 
= (-1)^{e (p_1-1)(p_2-1)/4} \chi^{\alpha_0 - N (N-1)/2-p_1 (p_1-1)/2-p_2
(p_2-1)/2} \zeta_{24 N}^{-24 r e N ((N-1)/2 - (p_1+p_2-2)/2)}.$$
Since $p_1+p_2-2$ is even, we obtain
$$\vartheta
= \zeta_{24}^{3 (24 r e)}
\chi^{\alpha_0 - N (N-1)/2-p_1 (p_1-1)/2-p_2 (p_2-1)/2
+ (N-1)/2 - (p_1+p_2-2)/2}
$$
$$=\chi^{\alpha_0-3 - N (N-1)/2-p_1 (p_1-1)/2-p_2 (p_2-1)/2
+ (N-1)/2 - (p_1+p_2-2)/2}$$
$$=\chi^{- N (N-1)/2-p_1 (p_1-1)/2-p_2 (p_2-1)/2 + (N-1)/2 - (p_1+p_2-2)/2},$$
and by inspection, this is always $1$.

\medskip
\noindent
(d) When $p_1 = p_2 = p$, we get
$$\prod_{\nu=1}^{p-1} C_{\nu}' = \prod_{\nu=1}^{p-1}
\chi^{\kappa(\nu)} p^{e/2} \varepsilon(\nu)^e \zeta_{24}^{e \theta(\nu)}
\zeta_{24p}^{-e \nu} q^{e (p-1)/12} (1 + \cdots)$$
$$= \chi^{p (p-1)/2} p^{e (p-1)/2}
\left(\prod_{\nu=1}^{p-1}\varepsilon(\nu)\right)^{e}
\zeta_{24}^{e\sum_{\nu=1}^{p-1} \theta(\nu)} \zeta_{24p}^{-e p
(p-1)/2} q^{2er} (1 + \cdots).$$
The quantity $\prod_{\nu=1}^{p-1}\varepsilon(\nu)$ is 1 for $p=2$;
when $p$ is odd
$$\prod_{\nu=1}^{p-1}\varepsilon(\nu) = \legendre{(-1)^{p-1}
(p-1)!}{p} = \legendre{-1}{p}$$
using Wilson's theorem.

When $p = 2$, $e = 8$, $\delta = 3$, we find
\begin{eqnarray*}
\vartheta &= &
\chi^{\alpha_0 - N (N-1)/2} \zeta_{24N}^{-24 r e N(N-1)/2}
\chi^{p (p-1)/2} p^{e (p-1)/2}
\zeta_{24}^{e\sum_{\nu=1}^{p-1} \theta(\nu)} \zeta_{24p}^{-e p (p-1)/2}\\
&=&\chi^{0 - 6} \zeta_{96}^{-8 \cdot 4 \cdot 3 / 2}
\chi^{1} 2^{4} \zeta_{48}^{-8} = 2^4.
\end{eqnarray*}
When $p$ is odd
$$\vartheta = 
\legendre{-1}{p}^e \chi^{\alpha_0 - N (N-1)/2 + p (p-1)/2} p^{e (p-1)/2}
\zeta_{24N}^{-24 r e N(N-1)/2}
\zeta_{24}^{e (-(p-1)/2 + \sum_{\nu=1}^{p-1} \theta(\nu))}.$$
Now:
$$\sum_{\nu=1}^{p-1} \theta(\nu) = \sum_{\nu=1}^{p-1} p \nu (1-\mu^2)
+ (-3 p + 2 + v) \mu - 3 + 3 p =
(p-3p+2) S_0 + 3 (p-1)^2 + \sum_{\nu=1}^{p-1} -p \nu \mu^2 + v \mu$$
where $S_0 = \sum_{\nu=1}^{p-1} \nu$. Using $1 = -\mu \nu + v p$, the
sum becomes
$$\sum_{\nu=1}^{p-1} p \mu (1 - v p) + v \mu
= \sum_{\nu=1}^{p-1} \mu (p - v p^2 + v) \equiv p S_0\bmod 24$$
in all cases: when $p > 3$, $p^2 \equiv 1 \bmod 24$; when $p=3$, 
$\nu = 1$ (resp. $\nu = -1$) leads to $\mu = -1$, $v = 0$
(resp. $\mu = 1$, $v = 0$). Therefore, the exponent of
$\zeta_{24}$ is 
$$\equiv e (-(p-1)/2 + (-p+2) S_0 + 3 (p-1)^2) \equiv -e
(p-7)(p-1)^2/2 \bmod 24,$$
so that
$$\vartheta = 
\legendre{-1}{p}^e p^{e (p-1)/2}
\zeta_{24N}^{-24 r e N(N-1)/2}
\chi^{\alpha_0 - N (N-1)/2 + p (p-1)/2}
\zeta_{24}^{-e (p-1)^2 (p-7)/2}
$$
$$= \legendre{-1}{p}^e p^{e (p-1)/2}
\chi^{\alpha_0 + (N-1)/2 - N (N-1)/2 + p (p-1)/2 + (p-7)/2}
= \legendre{-1}{p}^e p^{e (p-1)/2} \chi^{\alpha_0 - (p^4 - 3 p^2 - 8)/2}.$$
For instance, when $p=3$, $e=3$, $\delta = 2$, we find
$\vartheta = (-1) 3^3 (-1)^{1 - 31} = -3^3$.
More generally, as soon as $p > 3$,
$$\vartheta = \legendre{-1}{p}^e p^{e (p-1)/2} \chi^{\alpha_0 - 3} =
\legendre{-1}{p}^e p^{e (p-1)/2}$$
since $p^2 \equiv 1\bmod 24$ and the fact already used that
$\chi^{\alpha_0} = \chi^{-3}$.
$\Box$

\subsection{Computing the modular equations using series expansions}

There a variety of methods to compute the modular equations. For large
computations, it is possible to use suitably modified versions of
\cite{Enge09} or \cite{BrLaSu10}. Also, we can use resultants in the
same spirit as in the remark at the end of Section \ref{sct-weber},
noting that $\w_{p_1, p_2}^s = (\w_{p_1}(z)/\w_{p_1}(z/p_2))^s$.

Here, we content ourselves to use series expansions and nice formulas
that can help us for small cases. Also, this will add new properties
to our equations.

Looking carefully at the expression for $\mathcal{S}_k$, we see that
the terms in $C_1$, $C_2$ or $C$ cannot contribute to the modular
equation, since they have positive order. Therefore, we need only
consider the expansions of ${B'}^k$ and $\mathcal{S}_{A, k}$. Doing
this, we see that the useful terms for $\mathcal{S}_{A, k}$ are for $j
\leq -k t N'/N$. Since ${B'}^k = q^{-rek} (1 + \cdots)$ and $1\leq k \leq
\psi(N)$. We need at least $r e \psi(N)$ terms in the last coefficient.
Since $B'$ is the product and quotient of very sparse series, it might
be worthwhile to compute its powers by successive applications of
special routines handling this kind of computations. It is possible to
compute nice formulas for the $\mathcal{S}_{A, k}$, in the spirit of
the ones to come, but we do not need them.

A second algorithm consists in grouping
$$\Phi(F) = P_B(F) P_A(F) P_{C_1}(F) P_{C_2}(F)$$
and to compute $P_A$ (resp. $P_{C_1}$ and $P_{C_2}$) via its power
sums that are given in the preceding propositions.

Inspired by the approach of section \ref{ssct-wp}, 
the third algorithm uses the reciprocal polynomial, whose powers sums
will depend on the $C_{1, \nu}'$ and $C_{2, \nu}'$ only:
$$\Sigma_k = \sum_{\nu=0}^{p_2-1} \frac{1}{{C_{2, \nu}'}^k}
+ \sum_{\nu=0}^{p_1-1} \frac{1}{{C_{1, \nu}'}^k},$$
which is a process involving $p_1+p_2$. We will prove two useful
results (propositions \ref{SC2k} and \ref{SCk} below) to help us
compute these quantities.

\begin{proposition}\label{SC2k}
For all integers $k \neq 0$,
$$\mathcal{S}_{C_2, k} = p_2 \varepsilon_2^{k e} 
q^{k t/\delta}
\sum_{j \geq k t p_2'/p_2} c_{k, j p_2 - k t p_2'} q^j,$$
where $(p_2+1)/\delta = p_2'$ and the $c_{k, i}$ are explicitly given
in the proof.
\end{proposition}

\noindent
{\em Proof:} put $w = q^{1/p_2}$, $\zeta = \zeta_{p_2}$ and write
$$\frac{
	\eta\left(p_1 z\right)
	\eta\left(\frac{z}{p_2}\right)
  }{
	\eta(z)
	\eta\left(\frac{p_1 z}{p_2}\right)
  }
= 
\frac{
  \left(w^{p_1 p_2 /24} (1 + \sum_{i=1}^{\infty} a_i w^{p_1 p_2 i})\right)
  \left(w^{1/24} (1+\sum_{i=1}^{\infty} a_i w^{i})\right)
}
{
  w^{p_2 /24} \left(1 + \sum_{i=1}^{\infty} a_i w^{p_2 i}\right)
  w^{p_1 /24} \left(1 + \sum_{i=1}^{\infty} a_i w^{p_1 i}\right)
}
= w^r \mathcal{C}_{12}(w)
$$
with $\mathcal{C}_{12}(q) = 1 + \cdots \in \Z[[q]]$ (which is
  symmetrical in $p_1$ and $p_2$), from which
$$C_{2, \nu}' = \chi^{\theta_2(\nu)-\nu} 
\chi^{-\theta_2(\nu)} \varepsilon_2^e (w\zeta^\nu)^{r e}
  \mathcal{C}_{12}(w \zeta^\nu)^e.$$
and
$$\mathcal{S}_{C_2, k} = \varepsilon_2^{k e} w^{k r e} \sum_{\nu=0}^{p_2-1}
\chi^{-k \nu} \zeta^{k r e\nu} \mathcal{C}_{12}(w\zeta^\nu)^{e k}.$$
Writing $\mathcal{C}_{12}(w)^{ek} = \sum_{i=0}^\infty c_{k, i} w^i$
(remark this is valid irrespective of the sign of $k$),
the inner sum becomes
$$\sum_{i=0}^\infty c_{k, i} w^i
\sum_{\nu=0}^{p_2-1} (\chi^{-k} \zeta^{k r e + i})^\nu,$$
in which the root of unity is $\chi^{-k} \zeta^{k r e} =
 \zeta_{24}^{24 k r e} \zeta^{k r e} = \zeta_{24 p_2}^{24 k
 r e p_2 + 24 k r e} = \zeta_{24 p_2}^{24 k r e (p_2+1)}$. Now, we use
 the fact that $p_2 \equiv -1\bmod \delta$, so that $r e (p_2+1) = t
 p_2'$ where $p_2' = (p_2+1)/\delta$. The above sum is now
$$\sum_{i=0}^\infty c_{k, i} w^i \sum_{\nu=0}^{p_2-1} 
\left(\zeta^{k t p_2' + i}\right)^\nu
= p_2 \sum_{i\equiv -k t p_2'\bmod p_2}^\infty c_{k, i} w^i
= p_2 w^{-k t p_2'} \sum_{j \geq k t p_2'/p_2} c_{k, j p_2 - k t p_2'} q^j.$$
leading to the result. $\Box$

\begin{proposition}\label{SCk}
In case $p_1 = p_2 = p$, for all $k \neq 0$,
$$\mathcal{S}_{C, k} \in \left(q^{-1/24} \frac{\eta(p
z)^2}{\eta(z)}\right)^{e k} \Z[[q]] = q^{K k/\delta} \Z[[q]],$$
where all series are explicited in the proof.
\end{proposition}

\medskip
\noindent
{\em Proof:} One uses $\zeta = \zeta_p$ in
$$\mathcal{S}_{C, k} = p^{ek/2}
\;\left(\frac{\eta(p z)^2}{\eta(z)}\right)^{ek}\;
\sum_{\nu=1}^{p-1} \varepsilon(\nu)^{ek} \frac{\chi^{k \kappa(\nu)}
\zeta_{24}^{e k \theta(\nu)}}{\eta(z+\nu/p)^{ek}}
$$
$$=p^{ek/2}\;\left(\frac{\eta(p z)^2}{\eta(z)}\right)^{ek}\;
\sum_{\nu=1}^{p-1} \varepsilon(\nu)^{ek}
q^{-ek/24} \zeta^{-ek\nu/24}
\frac{\chi^{k \kappa(\nu)}\zeta_{24}^{e k \theta(\nu)}}{(1+\sum_{i=1}^\infty a_i q^i
\zeta^{i \nu})^{ek}}$$
$$=p^{ek/2} \;\left(q^{-1/24} \frac{\eta(p z)^2}{\eta(z)}\right)^{ek}\;
\sum_{\nu=1}^{p-1} \varepsilon(\nu)^{ek}
\chi^{k \kappa(\nu)} \zeta_{24}^{e k \theta(\nu)}
\zeta^{-ek\nu/24} \mathcal{C}(q\zeta^\nu)^{ek}$$
where
$$\mathcal{C}(q) = \frac{1}{1+\sum_{i=1}^\infty a_i q^i}.$$
Writing $\mathcal{C}(q)^{ek} = \sum_{i=0}^\infty c_{k, i} q^i$ (same
remark on the sign of $k$),
the inner sum of the preceding relation is now
\begin{equation}\label{eq1}
\sum_{\nu=1}^{p-1} \varepsilon(\nu)^{ek}
\chi^{k \kappa(\nu)} \zeta_{24}^{e k \theta(\nu)}
\zeta^{-ek\nu/24}
\sum_{i=0}^\infty c_{k, i} (q \zeta^\nu)^i = \sum_{i=0}^\infty c_{k,
i} q^i \sum_{\nu=1}^{p-1} \varepsilon(\nu)^{ek}
\chi^{k \kappa(\nu)} \zeta_{24}^{e k \theta(\nu)} (\zeta^{-ek/24}
\zeta^i)^\nu.
\end{equation}
Let's treat the case $p=2$ first, with $e=8$. We get
$$\mathcal{S}_{C, k}=2^{4 k}
\;\left(\frac{q^{-1/24} \eta(2 z)^2}{\eta(z)}\right)^{8 k}\;
\sum_{i=0}^\infty c_{k, i} q^i (\zeta_2^{-k} \zeta_2^i)
= (-2^4)^k \;\left(q^{-1/24} \frac{\eta(2 z)^2}{\eta(z)}\right)^{8k}\;
\sum_{i=0}^\infty c_{k, i} (-q)^i.$$

For $p$ odd, the root of unity in the inner sum of (\ref{eq1}) is
$$\varepsilon(\nu)^{ek}
\zeta_{24p}^{e k (p (-(p-1)^2 \kappa(\nu) + \theta(\nu)) - \nu)}
(\zeta^i)^\nu,$$
the exponent of $\zeta_{24p}$ being
$$p (-(p-1)^2 \mu + p\nu\,(1-{\mu}^{2})+ \left( -3\,p+2+v \right)
\mu-3+3\,p) - \nu.$$
When $p=3$ and $e=3$, we find
$$\zeta_{24}^{k (-32\,\mu+8\,\nu-8\,\nu\,{\mu}^{2}+18)} 
= (\zeta_{12}^{-16\,\mu+4\,\nu-4\,\nu\,{\mu}^{2}+9})^k
= (\zeta_4^3 \cdot \zeta_3^\nu)^k,$$
leading to
$$\mathcal{S}_{C, k}=
3^{3k/2} \zeta_4^{3k} 
\;\left(q^{-1/24} \frac{\eta(3 z)^2}{\eta(z)}\right)^{3k}\; 
\sum_{i=0}^\infty c_{k, i} q^i \sum_{\nu=1}^{2} \varepsilon(\nu)^{k}
(\zeta_3^{i+k})^\nu.$$
When $k$ is even, this boils down to
$$\mathcal{S}_{C, k}= (-3)^{3k/2}
\;\left(q^{-1/24} \frac{\eta(3 z)^2}{\eta(z)}\right)^{3k}\; 
\left(2 \sum_{i\equiv -k \bmod 3} c_{k, i} q^i
-\sum_{i\not\equiv -k \bmod 3} c_{k, i} q^i\right)
$$
$$= (-3)^{3k/2}
\;\left(q^{-1/24} \frac{\eta(3 z)^2}{\eta(z)}\right)^{3k}\; 
\left(3 \sum_{i\equiv -k \bmod 3} c_{k, i} q^i
-\sum_{i=0}^\infty c_{k, i} q^i\right)
.$$
When $k$ is odd
$$\mathcal{S}_{C, k}=
3^{3k/2} \zeta_4^{3k}
\;\left(q^{-1/24} \frac{\eta(3 z)^2}{\eta(z)}\right)^{3k}\;
\sum_{i=0}^\infty c_{k, i} q^i \sum_{\nu=1}^{2} \varepsilon(\nu)^{k}
(\zeta_3^{i+k})^\nu$$
and 
$$\sum_{\nu=1}^{2} \varepsilon(\nu)^{k} (\zeta_3^{i+k})^\nu =
- \zeta_3^{i+k} + \zeta_3^{2(i+k)}
=\left\{\begin{array}{cl}
0 & \text{if } i+k\equiv 0\bmod 3,\\
(-1)^{(i+k) \bmod 3} \sqrt{-3} & \text{otherwise},\\
\end{array}\right.$$
from which
$$\mathcal{S}_{C, k}=
(-3)^{(3k+1)/2}
\;\left(q^{-1/24} \frac{\eta(3 z)^2}{\eta(z)}\right)^{3k}\;
\sum_{i+k\not\equiv 0\bmod 3}^\infty (-1)^{(i+k)\bmod 3} c_{k, i} q^i.$$

When $p>3$, we get
$$(\zeta_{24}^{-(p-1)^2 \mu +
p\nu\,(1-{\mu}^{2})+ \left( -3\,p+2+v \right)
\mu-3+3\,p})^{ek} (\zeta_{24p}^{-ek +24 i})^\nu
=\left(\zeta_{24}^{\left( \nu-\nu\,{\mu}^{2}-\mu+3
\right) p+\mu\,v-3}\right)^{ek} (\zeta_{24p}^{-ek +24 i})^\nu$$
using $p^2 \equiv 1 \bmod 24$. We simplify this as
$$\left(\zeta_{24}^{p (\nu+3) - 3}\right)^{ek}
(\zeta_{24p}^{-ek +24 i})^\nu
=\zeta_{8}^{e k (p-1)}
\left(\zeta_{24p}^{ek (p^2-1) + 24 i}\right)^\nu.$$
Write $p^2-1 = 24 p'$ to obtain
$\zeta_{8}^{e k (p-1)} (\zeta_p^{e k p' + i})^\nu$.

When $k$ is even, this gives
$$\mathcal{S}_{C, k}=p^{ek/2} \zeta_{8}^{e k (p-1)}
\;\left(q^{-1/24} \frac{\eta(p z)^2}{\eta(z)}\right)^{ek}\; \sum_{i=0}^\infty
c_{k, i} q^i \sum_{\nu=1}^{p-1} (\zeta_p^{ekp'+i})^\nu$$
$$=
p^{ek/2} \zeta_{8}^{e k (p-1)}
\;\left(q^{-1/24} \frac{\eta(p z)^2}{\eta(z)}\right)^{ek}\;
\left((p-1) \sum_{i+ekp'\equiv 0\bmod p}^\infty c_{k, i} q^{i} -
\sum_{i+ekp'\not\equiv 0\bmod p} c_{k, i} q^i\right)$$
$$=
p^{ek/2} \zeta_{8}^{e k (p-1)}
\;\left(q^{-1/24} \frac{\eta(p z)^2}{\eta(z)}\right)^{ek}\;
\left(p\; \sum_{i+ekp'\equiv 0\bmod p}^\infty c_{k, i} q^{i} -
\sum_{i=0}^\infty c_{k, i} q^i\right).$$

When $k$ is odd, remarking that $e$ is always odd from Table
\ref{thp1p2}, the sum is now
$$\mathcal{S}_{C, k}=p^{ek/2} \zeta_{8}^{e k (p-1)}
\;\left(q^{-1/24} \frac{\eta(p z)^2}{\eta(z)}\right)^{ek}\; \sum_{i=0}^\infty
c_{k, i} q^i \sum_{\nu=1}^{p-1} \varepsilon(\nu) (\zeta_p^{i + ekp'})^\nu.$$
But
$\sum_{\nu=1}^{p-1} \varepsilon(\nu) (\zeta_p^{i + ekp'})^\nu = 0$ when
$i+e k p' \equiv 0 \bmod p$ since there are the same number of quadratic
residues and non-quadratic residues modulo $p$. When $i + e k p'\not\equiv
0\bmod p$, $\zeta_p^{i + ek p'}$ is a primitive $p$-th root of
unity. Remember that \cite[Ch. 6]{IrRo82}
$$\sum_{x\text{ residue}} \zeta_p^x -
\sum_{x\text{ non residue}} \zeta_p^x = \sqrt{\legendre{-1}{p}
p}.$$
Let $g$ be a generator of $\ZNS{p}$. If $u$ is an integer, then
$$\sum_{x\text{ residue}} (\zeta_p^{g^u})^x -
\sum_{x\text{ non residue}} (\zeta_p^{g^u})^x =
(-1)^u \sqrt{\legendre{-1}{p} p}.$$
When $i + e k p'\not\equiv 0\bmod p$, we set $\Omega(i+ekp') = u$ such
that $g^u \equiv i + e k p' \bmod p$. Then
$$\mathcal{S}_{C, k}=
\legendre{-1}{p}\sqrt{\legendre{-1}{p}}\zeta_{8}^{e k (p-1)}p^{(ek+1)/2} 
\;\left(q^{-1/24}\frac{\eta(p z)^2}{\eta(z)}\right)^{ek}\; 
\sum_{i=0, i+ek p' \not\equiv 0\bmod p}^\infty
(-1)^{\Omega(i+ekp')}c_{k, i} q^i.$$
When $p\equiv 1\bmod 4$, the first terms simplify to $\zeta_2^{e k
(p-1)/4} = (-1)^{(p-1)/4}$; when $p\equiv 3\bmod 4$, we get
$-\zeta_8^{2 + e k (p-1)} = -\zeta_4^{1 + e k (p-1)/2} = (-1)^{(3 + e k
(p-1)/2)/2}$.

As a last point, the dominant term of $\mathcal{S}_{C, k}$ is $q^{k e
(p-1)/12}$. When $p=2$ and $e=8$, this is $2k/3$, whereas $r e =
1/3$; when $p = 3$, $e = 3$, we get $k/2$, whereas $r e = 1/2$. For $p
> 3$, we have $e=1$ and we compare $(p-1)/12$ and $r e = (p-1)/2 \cdot
(p-1)/12$. Looking at the valuation of $2$ and $3$, we deduce that $r
e = t/\delta$ and $(p-1)/12 = p'/\delta$. $\Box$

\subsection{Tables of equations for double $\eta$-quotients}

$$\Phi[\w_{2, 2}^8](F, G_2)=
{F}^{6}-G_2\,{F}^{5}+208\,{F}^{3}+31\,G_2\,{F}^{2}+G_2^{2}F+16,$$
$$\Phi[\w_{3, 3}^3](F, G_3)=
{F}^{12}-G_3\,{F}^{11}-522\,{F}^{10}+27\,G_3\,{F}^{9}-10557\,{F}^{8}-162\,G_3\,{F}^{7}-14076\,{F}^{6}-18\,G_3\,{F}^{5}$$
$$-9801\,{F}^{4}+163\,G_3\,{F}^{3}+(486\,-G_3^2){F}^{2}-9\,G_3\,F-27.$$
$$\Phi[\w_{3, 7}](F, G_3)=
{F}^{32}-G_3\,{F}^{31}-514\,{F}^{30}+21\,G_3\,{F}^{29}-12585\,{F}^{28}-147\,G_3\,{F}^{27}-25158\,{F}^{26}+322\,G_3\,{F}^{25}$$
$$-5103\,{F}^{24}+378\,G_3\,{F}^{23}+80556\,{F}^{22}-1638\,G_3\,{F}^{21}-21994\,{F}^{20}-28136\,{F}^{18}+1620\,G_3\,{F}^{17}+25650\,{F}^{16}$$
$$-252\,G_3\,{F}^{15}-3944\,{F}^{14}-322\,G_3\,{F}^{13}-14938\,{F}^{12}+22\,G_3\,{F}^{11}-(G_3^{2}-2940)
{F}^{10}-10\,G_3\,{F}^{9}+1953\,{F}^{8}$$
$$+G_3\,{F}^{7}-462\,{F}^{6}+7\,G_3\,{F}^{5}+15\,{F}^{4}-G_3\,{F}^{3}-10\,{F}^{2}+1.$$

\section{Conclusion}

We have studied modular equations involving $\gamma_2$ and $\gamma_3$
for double $\eta$-quotients. As a result, more compact modular
equations can be stored and used, with application to the SEA
algorithm (see for instance \cite{Morain95a}), or CM computations, as
motivated for instance by \cite{RuSi10} (see \cite{Morain11}).

It seems natural to conjecture that more
general functions can exhibit the same properties. Experiments can be
conducted on Newman functions, using for instance the resultant
approach, leading to new instances of the theorems. This will
described in another article.

\iffalse
\bibliographystyle{plain}
\bibliography{morain}
\else
\def\noopsort#1{}\ifx\bibfrench\undefined\def\biling#1#2{#1}\else\def\biling#1%
#2{#2}\fi\def\Inpreparation{\biling{In preparation}{en
  pr{\'e}paration}}\def\Preprint{\biling{Preprint}{pr{\'e}version}}\def\Draft{%
\biling{Draft}{Manuscrit}}\def\Toappear{\biling{To appear}{\`A para\^\i
  tre}}\def\Inpress{\biling{In press}{Sous presse}}\def\Seealso{\biling{See
  also}{Voir {\'e}galement}}\def\Editor{\biling{Ed.}{R{\'e}d.}}

\fi

\end{document}